\documentclass[a4paper]{amsart}

\newtheorem{theorem}{Theorem}[section]
\newtheorem{lemma}[theorem]{Lemma}
\newtheorem{prop}[theorem]{Proposition}
\newtheorem{coroll}[theorem]{Corollary}

\theoremstyle{definition}
\newtheorem{definition}[theorem]{Definition}

\theoremstyle{remark}
\newtheorem{remark}[theorem]{Remark}

\numberwithin{equation}{section}

\newcommand{\f}[1][p]{\mathbb{F}_{\scriptstyle #1}}
\newcommand{\z}{\mathbb{Z}}
\newcommand{\C}{\mathbb{C}}
\newcommand{\proj}{\mathbb{P}}

\newcommand{\Res}{\mathop{\text{\rm Res}}\nolimits}
\newcommand{\Cor}{\mathop{\text{\rm Cor}}\nolimits}
\newcommand{\Inf}{\operatorname{Inf}}
\newcommand{\Bild}{\mathop{\text{\rm Im}}\nolimits}
\newcommand{\Hom}{\mathop{\text{\rm Hom}}\nolimits}

\newcommand{\coho}[3][*]{\mathrm{H}^{#1}(#2,#3)}
\newcommand{\hz}[2][*]{\coho[#1]{#2}{\z}}
\newcommand{\hf}[2][*]{\coho[#1]{#2}{\f}}

\newcommand{\Chop}{\mathit{Ch}}

\newcommand{\GL}{\mathit{GL}}

\newcommand{\Chbarop}{\overline{\Chop}}
\newcommand{\Chbar}[1]{\Chbarop (#1)}

\newcommand{\DrV}[2]{D_{#1}(#2)}
\newcommand{\Qnr}[2]{D_{#1,#2}}
\newcommand{\Muiop}{V}
\newcommand{\Mui}[2]{\Muiop(#1;#2)}
\newcommand{\MuiRel}[2]{V(#1,#2)}
\newcommand{\Hyp}[2]{\mathit{Hyp}(#1,#2)}
\newcommand{\Pp}{P}
\newcommand{\Pn}[1]{\Pp_{#1}}
\newcommand{\Ee}{E}
\newcommand{\En}[1]{\Ee_{#1}}
\newcommand{\vsp}[1]{\mathcal{#1}}
\newcommand{\Tt}{\vsp{T}}
\newcommand{\Uu}{\vsp{U}}
\newcommand{\Vv}{\vsp{V}}
\newcommand{\Ww}{\vsp{W}}
\newcommand{\vctr}[1]{#1}
\newcommand{\vv}{\vctr{v}}
\newcommand{\ww}{\vctr{w}}
\newcommand{\Iphi}{I_{\phi}}
\newcommand{\Ibeta}{I_{\beta_n}}
\newcommand{\Ss}{S}
\newcommand{\Rr}{R}
\newcommand{\s}{s}
\newcommand{\zni}[2][n]{z_{#1}^{(#2)}}
\newcommand{\ZH}{Z_H}
\newcommand{\zdi}[1]{{z'_1}^{(#1)}}

\begin{document}

\title{Transfer and Chern Classes for Extraspecial $p$-Groups}

\author[D.~J. Green]{David John Green}
\address{Inst.\@ f.\@ Exp.\@ Math.\@ \\ Ellernstr.~29 \\ D--45326 Essen \\
Germany}
\email{david@exp-math.uni-essen.de}
\thanks{The first author was supported by the Deutsche Forschungsgemeinschaft
Schwerpunktprogramm ``Algorithmische Zahlentheorie und Algebra''.}

\author[P.~A. Minh]{Pham Anh Minh}
\address{Fachbereich 6 Mathematik \\ Universit\"at Essen \\ D--45117 Essen \\
Germany}
\curraddr{Department of Mathematics \\ University of Hue \\ Dai hoc Tong
hop Hue \\ Hue \\ Vietnam}
\thanks{The second author held a DAAD fellowship.}

\date{3 March 1997}

\begin{abstract}
In the cohomology ring of an extraspecial $p$-group, the subring generated
by Chern classes and transfers is studied.  This subring is
strictly larger than the Chern subring, but still not the whole cohomology
ring, even modulo nilradical.  A formula is obtained relating Chern classes
to transfers.
\end{abstract}

\maketitle

\section*{Introduction}
Methods to determine the cohomology ring of a finite group almost
always presuppose that the cohomology of the Sylow $p$-subgroups is
known.  Calculating the cohomology ring of a $p$-group is however
a delicate and difficult task.  The extraspecial $p$-groups of exponent~$p$
are in some sense the minimal difficult cases: minimal because their proper
quotients are all elementary abelian, and
their automorphism groups are very large.
For this reason, many papers have been written, investigating their cohomology.
Developments up till 1991 are surveyed in the paper~\cite{BC92}.
In particular, M.~Tezuka and N.~Yagita obtained the prime ideal spectrum of the
cohomology ring.

The usual method to calculate the cohomology of a $p$-group is to write
the group as an extension, and solve the associated Lyndon--Hochschild--Serre
spectral sequence.  For the extraspecial $p$-groups however, such spectral
sequences are intractable, and one is forced to look for other techniques.
Now, standard constructions such as transfer (or corestriction) from subgroups
and taking Chern classes of group representations provide us with a large
number of cohomology classes.  So many in fact, that for any $p$-group the
classes provided by these two constructions generate a subring that has the
same prime ideal spectrum as the cohomology ring (see~\cite{VCh}).
In this paper, we study this subring in the case of the extraspecial $p$-groups,
and ask whether it is the whole cohomology ring.

Actually, these constructions yield very few odd-dimensional classes,
and so it is rather more realistic to ask whether we obtain the whole
cohomology ring modulo nilradical.  At least for mod-$p$ cohomology,
Proposition~\ref{prop:Chbar} answers this question too in the negative.
For integral cohomology however, the problem remains open, and the significance
of Corollary~\ref{cor:bigger} is that this subring is the biggest studied
to date in the cohomology ring of an extraspecial $p$-group.

This last result is proved using cohomology classes which we
denote~$\chi_{r,\phi}$.  They are constructed in a manner foreshadowed
in~\cite{Minh:transfer}.  Take a product of Chern classes for the group
of order~$p^{2n-1}$.  By inflation and then corestriction, obtain a cohomology
class~$\chi_{r,\phi}$ for the group of order~$p^{2n+1}$.
In Theorem~\ref{thm:5.2}, an elegant formula is obtained relating the Chern
classes and the~$\chi_{r,\phi}$, and in Theorem~\ref{thm:pthPower} we show
that the $p^{\text{th}}$ power of any Chern class or any~$\chi_{r,\phi}$
lies in the subring generated by top Chern classes: this is the subring
Tezuka and Yagita used to obtain the spectrum of the cohomology ring.


We are very grateful to Bruno Kahn for interesting discussions; and to
H\'el\`ene Esnault and Eckhart Viehweg, who arranged for the second
author to visit Essen.

\section{A relation between Dickson invariants}
We shall assume that the reader is familiar with the Dickson invariants:
such familiarity may be acquired by consulting
Benson's book~\cite{Benson:PolyInvts}, for example.  

Let~$\Vv$ be an $m$-dimension $\f$-vector space, and let $0 \leq r \leq m-1$.
We shall write~$\DrV{r}{\Vv}$ for the Dickson invariant in degree $p^m - p^r$
in~$S(\Vv)$, or just~$\Qnr{m}{r}$ if~$\Vv$ is clear from the context.
Recall that, for an indeterminate~$X$, we have the equation
\begin{equation}
\label{eqn:DicksonTheorem}
\Mui{\Vv}{X} = X^{p^m} + \sum_{r=0}^{m-1} (-1)^{m-r} \DrV{r}{\Vv} X^{p^r}
\quad ,
\end{equation}
where $\Mui{\Vv}{X}$~in $S(\Vv)\lbrack X \rbrack$ is defined by
\begin{equation}
\label{eqn:Muidef}
\Mui{\Vv}{X} = \prod_{\phi \in \Vv} (X - \phi) \quad .
\end{equation}
The polynomial~$\Mui{\Vv}{X}$ is called the Mui invariant, as it is the
most important of a family of invariants for subgroups of~$\GL_{m+1}$
studied by Mui~\cite{Mui75}.  Note that~$\Mui{\Vv}{X}$ is $\f$-linear as
a function of~$X$.
By convention, we define $\Qnr{m}{r} = 0$ if $r < 0$, and $\Qnr{m}{m} = 1$.

\begin{lemma}
\label{lemma:DicksonRelation}
Let~$\Vv$ be an $m$-dimensional $\f$-vector space, let~$\Uu$ be a proper
subspace of~$\Vv$, and let~$s$ be the codimension of $\Uu$~in $\Vv$.
Let $\Hyp{\Vv}{\Uu}$ denote the set of all hyperplanes in~$\Vv$ which
contain~$\Uu$: this set is clearly nonempty.  Then
\begin{equation}
\sum_{\Ww \in \Hyp{\Vv}{\Uu}} \Mui{\Ww}{X} = \Mui{\Uu}{X}^{p^{s-1}} \quad ,
\end{equation}
and it follows that
\begin{equation}
\sum_{\Ww \in \Hyp{\Vv}{\Uu}} \DrV{r-1}{\Ww} = \DrV{r-s}{\Uu}^{p^{s-1}}
\quad \text{for all $1 \leq r \leq m$.}
\end{equation}
\end{lemma}

\begin{proof}
Let~$\Tt$ be a complementary subspace of $\Uu$~in $\Vv$,
so $\dim(\Tt) = s$.  The elements of~$\Hyp{\Vv}{\Uu}$ are those
subspaces $\Ww$~of $\Vv$ such that $\Ww \cap \Tt$ is a hyperplane in~$\Tt$
and $\Ww = \Uu \oplus (\Ww \cap \Tt)$.  So choosing such a~$\Tt$ induces
a bijection between~$\Hyp{\Vv}{\Uu}$ and the projective space $\proj \Tt^*$ of
lines in the dual space of~$\Tt$.

Let $\tau_1, \ldots, \tau_s$ be a basis for~$\Tt$.
Observe that~$\Mui{\Uu}{\tau_1}$, \dots,~$\Mui{\Uu}{\tau_s}$,
all homogeneous of degree~$p^{n-s}$,
are algebraically independent over~$\f$:
for the image of $\Mui{\Uu}{\tau_i}$~in $S(\Vv/\Uu) \cong S(\Tt)$
is~$\tau_i^{p^{n-s}}$.

Denote by~$\underline{\lambda}$ an $s$-tuple $(\lambda_1,\ldots,\lambda_s)$~in
$\f^s$.  Given~$\phi \in \proj \Tt^*$, say that~$\underline{\lambda}$
belongs to~$\ker(\phi)$ if $\sum_{i=1}^s \lambda_i \tau_i$ does.
%
Using Eqn.~(\ref{eqn:Muidef}) and the linearity of $\Mui{\Uu}{{-}}$, we obtain
\begin{align}
\sum_{\Ww \in \Hyp{\Vv}{\Uu}} \Mui{\Ww}{X}
&= \sum_{\phi \in \proj \Tt^*} \prod_{t \in \ker (\phi)} \Mui{\Uu}{X - t} \\
&= \sum_{\phi} \prod_{\underline{\lambda} \in \ker(\phi)}
  \left( \Mui{\Uu}{X} - \sum_{i=1}^s \lambda_i \Mui{\Uu}{\tau_i} \right)
\label{eqn:bigsum}
\quad .
\end{align}
Consider the
expression~(\ref{eqn:bigsum}).
It is a polynomial
in~$\Mui{\Uu}{X}$, $\Mui{\Uu}{\tau_1}$, \dots,~$\Mui{\Uu}{\tau_s}$, and as such
is homogeneous of degree~$p^{s-1}$.  Treat~$\Mui{\Uu}{X}$ as the main
variable: then there is no constant term.  Moreover, the  coefficients are
polynomials in~$\Mui{\Uu}{\tau_1}$,
\dots,~$\Mui{\Uu}{\tau_s}$ which are invariant under the natural action
of $\GL_s(\f)$.  So the coefficients must be polynomials
in the Dickson invariants in these~$s$ variables: but there are no such
polynomials in positive degree less than~$p^{s-1}$.

Therefore the coefficient in~(\ref{eqn:bigsum}) of~$\Mui{\Uu}{X}^j$ is
zero except when $j = p^{s-1}$.  This coefficient is the size of~$\proj \Tt^*$,
congruent to~$1$ modulo~$p$.
\end{proof}

\section{Chern classes and extraspecial $p$-groups}
\label{section:Chern}
The integral (group) cohomology of the unitary group~$U(n)$
is a polynomial algebra with~$n$ generators.  A unitary representation of
a finite group~$G$ pulls these generators back to the integral cohomology
of~$G$; the images of the generators are the Chern classes of the
representation.  We refer the reader to~\cite{CBT} for more information about
Chern classes and their properties.

Let~$p$ be an odd prime.
For $n \geq 1$, denote by~$\Pp = \Pn{n}$ the extraspecial $p$-group of
order~$p^{2n+1}$ and exponent~$p$.  This fits into a central extension
\begin{equation}
1 \rightarrow Z \rightarrow \Pn{n} \rightarrow \En{n} \rightarrow 1 \quad ,
\end{equation}
where~$Z = Z(\Pp)$ is cyclic of order~$p$, and~$\Ee = \En{n}$ is
elementary abelian of $p$-rank~$2n$.  We may identify~$Z$ with~$\f$
and view~$\Ee$ as a $2n$-dimensional $\f$-vector space: then
the commutator map on~$\Pp$ induces a nondegenerate symplectic
bilinear form $\Ee \otimes_{\f} \Ee \rightarrow Z$.

Ignoring the~$p^{2n}$ linear characters, the remaining irreducible
characters of~$\Pp$ all have degree~$p^n$, and are distinguished by their
restrictions to~$Z$.  Pick a nontrivial linear character $\gamma$~of $Z$,
and define~$\rho_{\gamma}$ to be the unique irreducible representation
of~$\Pp$ whose restriction to~$Z$ has character~$p^n \gamma$.  The other
degree~$p^n$ representations are then $\psi^j (\rho_{\gamma})$ for
$2 \leq j \leq p-1$.

Since~$\Ee$ is the abelianization of~$\Pp$, Chern classes of one-dimensional
representations lie in the image of inflation from $\Ee$~to $\Pp$.  The
Chern classes of~$\psi^j(\rho_{\gamma})$ are scalar multiples of those
of~$\rho_{\gamma}$.

\medskip
Pick an embedding of $\f$~in $\C^{\times}$: for each elementary abelian
$p$-group~$A$, this allows us to identify $\Hom(A,\C^{\times})$ with
the dual space~$A^*$.

The maximal elementary abelian subgroups $M$~of $\Pp$ all have
$p$-rank $n + 1$, and are permuted transitively by the automorphism group
of~$\Pp$.  Each contains~$Z$, and the map $M \mapsto I = M / Z$
induces a bijection between the~$M$ and the maximal totally isotropic
subspaces $I$~of $\Ee$.  Note that the dual~$I^*$ is naturally a subspace
of~$M^*$.

View~$\gamma$ as an element of $\hz[2]{Z}$, since this isomorphic to
$\Hom(Z,\C^{\times})$.  The restriction of~$\rho_{\gamma}$ to~$M$ decomposes
as the direct sum of all representations whose restriction to~$Z$ has
character~$\gamma$.  The first Chern classes of these summands are all the
classes in $\hz[2]{M}$ whose restriction to~$Z$ is~$\gamma$; pick~$\bar\gamma$
to be one of these first Chern classes, and observe that the kernel
in $\hz[2]{M} \cong M^*$ of restriction to~$Z$ is~$I^*$.  From the
definition of~$\Muiop$, the total Chern class
of~$\rho_{\gamma}$ therefore restricts to~$M$ as $\Mui{I^*}{1 + \bar\gamma}$,
independent of the choice of~$\bar\gamma$.
Define $\zeta_n \in \hz[2 p^n]{\Pp}$ and, for $0 \leq r \leq n - 1$, define
$\kappa_r = \kappa_{n,r}$ in $\hz[2(p^n - p^r)]{\Pn{n}}$ by
\begin{equation}
\zeta_n = c_{p^n}(\rho_{\gamma}) \qquad 
\kappa_{n,r} = (-1)^{n-r} c_{p^n - p^r} (\rho_{\gamma}) \quad .
\end{equation}
Then these are the only non-nilpotent Chern classes of~$\rho_{\gamma}$, and
\begin{equation}
\Res^{\Pp}_M \zeta_n = \Mui{I^*}{\bar\gamma} \qquad
\Res^{\Pp}_M (\kappa_r) = \DrV{r}{I^*} \quad .
\end{equation}

It is useful, and sensible, to define $\Pn{0}$ to be~$Z$ and
$\zeta_0$~to be~$\gamma$.  Also to define~$\kappa_{n,n}$ to be
$1 \in \hz[0]{\Pp}$.

\section{The new classes}
Since every group homomorphism $\Pp \rightarrow \f$ factors through~$\Ee$,
we may identify $\Ee^*$~with the $\f$-vector space $\Hom(\Pp, \f)$.
Each maximal subgroup $H \leq \Pp$ induces a one-dimensional subspace
of $\Hom(\Pp, \f)$, namely the subspace generated by any~$\phi$
with $\ker (\phi) = H$.  Consequently there is a natural bijection
between the set of maximal subgroups of~$\Pp$ and the projective
space~$\proj \Ee^*$.

We are interested in the corestriction map from maximal subgroups to~$\Pp$.
Let $\phi \in \proj \Ee^*$, and let $H = \ker(\phi)$ be the corresponding
maximal subgroup of~$\Pp$.  Then $H / Z$ is a $(2n - 1)$-dimensional
$\f$-vector space, carrying a symplectic bilinear form with one-dimensional
kernel.  The centre $\ZH$~of $H$ is elementary abelian of $p$-rank~$2$,
and $\ZH / Z$ is the kernel of the form on~$H / Z$.  The
maximal elementary abelian subgroups of~$H$ all have $p$-rank $n+1$ and all
contain~$\ZH$.  They are permuted transitively by the automorphism
group of~$H$.


\begin{definition}
\label{definition:chi}
For $n \geq 1$, for $0 \leq r \leq n-1$ and for $\phi \in \proj \En{n}^*$,
the class~$\chi_{r,\phi}$ in $\hz{\Pn{n}}$ with degree $2 (p^n - p^r)$ is
defined as follows.  Set $H = \ker(\phi)$ and pick a rank one
subgroup $A \not = Z$~of $\ZH$.  Choosing such an~$A$ induces
a split epimorphism $H \rightarrow H/A \cong \Pn{n-1}$.
Set
\begin{equation}
\chi_{r,\phi} =
  \Cor^{\Pn{n}}_H \Inf^H_{\Pn{n-1}} (\kappa_{n-1,r} \zeta_{n-1}^{p-1} )
  \quad .
\end{equation}
\end{definition}


\begin{lemma}
The class~$\chi_{r,\phi}$ is well-defined.  That is, it does not depend on
the choice of $A \leq \ZH$.
\end{lemma}

\begin{proof}
The inner automorphisms of~$\Pn{n}$ permute transitively all such
subgroups $A \leq \ZH$, always fixing~$Z$ pointwise and therefore sending
the~$\rho_{\gamma}$ for one~$\Pn{n-1}$ to the~$\rho_{\gamma}$ for the other.
Since~$H$ is normal in~$\Pn{n}$ and corestriction
commutes with conjugation, the result follows.
\end{proof}

\section{Restriction and the new classes}
In this section, we study restrictions of the classes~$\chi_{r,\phi}$.
We start however with a preparatory lemma.

\begin{lemma}
\label{lemma:MuiRel}
Let~$\Vv$ be an $\f$-vector space, and let~$\Uu$ be a hyperplane in~$\Vv$.
Pick~$\vv$ in $\Vv \setminus \Uu$.  Then the element
$\Mui{\Uu}{\vv}^{p-1}$ of~$S(\Vv)$ is invariant under all transformations
of~$\Vv$ which preserve~$\Uu$.  In particular,
it is independent of the choice of~$\vv$.
\end{lemma}

\begin{proof}
Any transformation which preserves~$\Uu$ acts on the coset space~$\Vv / \Uu$
as multiplication by some scalar in~$\f^{\times}$.  But, from
its definition, $\Mui{\Uu}{\vv}$ is invariant under all transformations
which preserve $\Uu$~and the coset $\vv + \Uu$.
\end{proof}

\begin{definition}
In the situation of Lemma~\ref{lemma:MuiRel}, we denote
by~$\MuiRel{\Uu}{\Vv}$
the invariant $\Mui{\Uu}{\vv}^{p-1}$.
\end{definition}

\begin{lemma}
\label{lemma:MuiRel2}
Let~$\Vv$ be an $\f$-vector space, and let~$\Uu$ be a subspace with
codimension two in~$\Vv$.  Then
\begin{equation}
\sum_{\Ww} \MuiRel{\Uu}{\Ww} = 0 \quad ,
\end{equation}
where the sum is over all hyperplanes $\Ww$~in $\Vv$ which themselves
contain~$\Uu$ as a hyperplane.
\end{lemma}

\begin{proof}
Let $\vv,\ww$ be a basis for a subspace of~$\Vv$ complementary to~$\Uu$.
Then
\begin{align}
\sum_{\Ww} \MuiRel{\Uu}{\Ww}
&= \sum_{\lbrack \lambda \colon \mu \rbrack \in \f \proj}
   \Mui{\Uu}{\lambda \vv + \mu \ww}^{p-1} \\
&= \sum_{\lbrack \lambda \colon \mu \rbrack} \left( \lambda \Mui{\Uu}{\vv}
   + \mu \Mui{\Uu}{\ww} \right)^{p-1} \quad .
\end{align}
This is an invariant of $\GL_2(\f)$ acting on the rank~$2$ polynomial
algebra generated by $\Mui{\Uu}{\vv}$~and $\Mui{\Uu}{\ww}$.  Its degree in
the generators is however~$p-1$, and there are no invariants other than zero
in this degree.
\end{proof}

We now investigate the image of~$\chi_{r,\phi}$ under restriction to each
maximal elementary abelian $M \leq \Pp$.  If~$M$ is contained
in $H = \ker (\phi)$, then $\ZH$~is contained in~$M$.  Define~$\Iphi$
to be the quotient $M / \ZH$.  Then~$\Iphi^*$ is a hyperplane in~$I^*$,
which is itself a hyperplane in~$M^*$.  Concretely, $\Iphi^*$~is 
the annihilator in $I^*$~of $\ZH / Z$ and also the annihilator
in $M^*$~of $\ZH$.

\begin{prop}
\label{prop:ResChi}
With the above notation, we have
\begin{equation}
\Res_M (\chi_{r,\phi}) = \begin{cases}
- \DrV{r}{\Iphi^*} \MuiRel{\Iphi^*}{I^*} & \text{if $M \leq \ker (\phi)$} \\
0 & \text{otherwise.}
\end{cases}
\end{equation}
\end{prop}

\begin{proof}
Since~$\chi_{r,\phi}$ was defined as a corestriction, we use the Mackey
formula to determine its restriction to~$M$.  Both $M$~and $H$ are normal
in~$\Pp$, and corestriction to~$M$ from any proper subgroup is the zero
map, at least in positive degree.  This proves the result when~$M$ is not
contained in~$H$.  So we may now assume that~$M \leq H$.

As in the definition of~$\chi_{r,\phi}$, choose a cyclic subgroup
$A$~of $\ZH$ such that $\ZH = A \times Z$.  Then~$M/A$ is a
maximal elementary abelian subgroup of $H/A \cong \Pn{n-1}$.
Note that $(M/A)^*$ is
a hyperplane in~$M^*$, and itself contains~$\Iphi^*$ as a hyperplane: namely,
the annihilator of~$Z$.
Then
\begin{align}
\Res_M (\chi_{r,\phi})
&= \sum_{g \in \Pn{n} / H} g^* \Res^H_M \Inf^H_{\Pn{n-1}}
   \kappa_{n-1,r} \zeta_{n-1}^{p-1} \\
&= \sum_{g \in \Pn{n} / H} \DrV{r}{\Iphi^*} 
   \MuiRel{\Iphi^*}{g^* (M/A)^*} \quad . \label{eqn:MA}
\end{align}
Now, there were~$p$ possible choices for~$A$, each of which yields a
different $(M/A)^*$.  The possible~$(M/A)^*$ are exactly those hyperplanes
in~$M^*$ which contain~$\Iphi^*$ but are not equal to~$I^*$. There are
permuted faithfully and transitively by $\Pn{n} / H$.  The result therefore
follows from Eqn.~(\ref{eqn:MA}) by Lemma~\ref{lemma:MuiRel2}.
\end{proof}

\section{Describing Chern classes in terms of corestrictions}
The Chern classes~$\kappa_{n,r}$ restrict to~$Z$ as zero.
A theorem of Carlson~\cite[\S10.2]{Evens:book} says that any such class
has some power which is a sum of corestrictions from proper subgroups.
In this section we shall derive a formula for~$\kappa_{n,r}$ in terms of
the corestrictions~$\chi_{r,\phi}$ and the image of inflation from~$\Ee$.
First however, we recall a well-known fact about Dickson invariants.

\begin{lemma}
\label{lemma:DicksonInduction}
Let~$\Vv$ be an $\f$-vector space, and $\Uu$~a hyperplane in~$\Vv$.  Then
\begin{equation}
\Mui{\Vv}{X} = \Mui{\Uu}{X}^p - \Mui{\Uu}{X} \MuiRel{\Uu}{\Vv} \quad ,
\end{equation}
and so, for $0 \leq r \leq n-1$,
\begin{equation}
\DrV{r}{\Vv} = \DrV{r-1}{\Uu}^p + \DrV{r}{\Uu} \MuiRel{\Uu}{\Vv} \quad .
\end{equation}
Here, recall that $\DrV{-1}{\Uu}$ is zero.
\qed
\end{lemma}

Let~$M_0$ be a maximal elementary abelian subgroup of~$\Pp$, and denote
by~$J_0$ the quotient $\Pp / M_0$, itself elementary abelian of rank~$n$.
Then~$J_0^*$ is the annihilator in $\Ee^*$~of $M_0$.  In addition,
we may view~$J_0^*$ as a subspace of $\hz[2]{\Ee}$.

\begin{theorem}
\label{thm:5.2}
For all $0 \leq r \leq n - 1$, the degree $2(p^n - p^r)$ class~$\nu_{r,M_0}$
in $\hz{\Pp}$ defined by
\begin{equation}
\nu_{r,M_0} = \kappa_{n,r} - \Inf^{\Pp}_{\Ee} \DrV{r}{J_0^*}
 + \sum_{\phi \in \proj J_0^*} \chi_{r,\phi}
\end{equation}
is nilpotent.
\end{theorem}

\begin{proof}
By a theorem of Quillen~\cite[Cor.~8.3.4]{Evens:book}, a class in the mod-$p$
cohomology ring of a finite group is nilpotent if and only if the
restriction to every elementary abelian $p$-subgroup is.  Now,
a class in integral cohomology reduces to zero in mod-$p$ cohomology only
if it is in the ideal generated by~$p$.  But in positive degree, all
such classes are nilpotent.  Hence it suffices to prove that the
restriction to every elementary abelian subgroup is zero.

Let~$M$ be a maximal elementary abelian subgroup of~$\Pp$.  As before,
define $I = M / Z$, and $\Iphi = M / (Z (\ker \phi))$ for $\phi \in \proj \Ee^*$.
Define~$\Ss$ to be the annihilator in $J_0^*$~of $M$, recalling that~$J_0^*$
may be viewed as a space of homomorphisms $\Pp \rightarrow \f$.  Then
$\Ss = \lbrace \phi \in \Ee^* \mid M,M_0 \in \ker (\phi) \rbrace$.

If $\Ss = \lbrace 0 \rbrace$, then
$\Res_M (\chi_{r,\phi}) = 0$ for every $\phi \in J_0^*$,
and both $\kappa_{n,r}$~and $\Inf^{\Pp}_{\Ee} \DrV{r}{J_0^*}$ restrict
to $M$~as $\DrV{r}{I^*}$.  So $\Res_M (\nu_{r,M_0}) = 0$ as required.

We may therefore assume that the $\f$-vector space~$\Ss$ has positive
dimension~$\s$, which implies that the number of elements of the projective
space $\proj \Ss$ is congruent to one modulo~$p$.
Denote by~$\Rr$ the subspace $\Res_M J_0^*$ of~$I^*$, and observe that~$\Ss$
is the kernel of the quotient map $J_0^* \rightarrow \Rr$.
Now, $\Ss$~contains by
Proposition~\ref{prop:ResChi} the set of those $\phi \in J_0^*$ such that
$\Res_M (\chi_{r,\phi})$ is nonzero.  Hence
\begingroup
\allowdisplaybreaks
\begin{align}
\Res_M \left( \displaystyle\sum_{\phi \in \proj J_0^*}
  \chi_{r,\phi} \right)%
& = - \sum_{\phi \in \proj \Ss} \DrV{r}{\Iphi^*} \MuiRel{\Iphi^*}{I^*}
\quad \text{by Proposition~\ref{prop:ResChi}}
\\
& = \sum_{\phi \in \proj \Ss}
        \bigl( \DrV{r-1}{\Iphi^*}^p - \DrV{r}{I^*} \bigr)
\quad \text{by Lemma~\ref{lemma:DicksonInduction}} \\
& = \DrV{r - \s}{\Rr}^{p^{\s}} - \DrV{r}{I^*}
\quad \text{by Lemma~\ref{lemma:DicksonRelation}} \\
& = \Res_M \bigl( \Inf^{\Pp}_{\Ee} \DrV{r}{J_0^*} - \kappa_{n,r} \bigr)
\quad .
\end{align}
\endgroup
Therefore $\Res_M (\nu_{r,M_0}) = 0$, as desired.
\end{proof}

\section{Taking $p^{\text{th}}$ powers of Chern classes}
In the next section we shall show that the $p^{\text{th}}$ powers of the
new classes~$\chi_{r,\phi}$ lie in the image of inflation from~$\Ee$,
at least modulo nilradical.  In preparation for this, we shall here prove
the same result for the $p^{\text{th}}$ powers of the Chern
classes~$\kappa_{n,r}$.

Recall that~$\Ee$ carries a nondegenerate symplectic form, say~$(,)$.
Pick a symplectic basis $A_1, \ldots, A_n$, $B_1, \ldots, B_n$ for~$\Ee$:
so $A_i \perp A_j$, $B_i \perp B_j$ and $(A_i, B_j) = \delta_{ij}$.
Let $\alpha_1, \ldots, \alpha_n$, $\beta_1, \ldots, \beta_n$ be the
corresponding dual basis for~$\Ee^*$, which we recall may be identified with
$\hz[2]{E}$.

\begin{definition}
For $r \geq 1$, define~$\zni{r}$ in $S(\Ee^*)$ by
\begin{equation}
\zni{i} = \sum_{i=1}^n ( \alpha_i \beta_i^{p^r} - \alpha_i^{p^r} \beta_i)
\quad .
\end{equation}
\end{definition}

By the work of Carlisle and Kropholler on symplectic invariants
(see \cite[\S8.3]{Benson:PolyInvts}), the~$\zni{i}$ are invariant
under symplectic transformations of~$\Ee$.  In particular, this means that
they do not depend upon the choice of symplectic basis.
Tezuka and Yagita proved that, at least for $p$~odd,
$\zni{1}, \ldots, \zni{n}$ is a regular sequence in~$S(E^*)$; and
that the ideal they generate contains every~$\zni{i}$ and is the
intersection with~$S(E^*)$ of the kernel of inflation
from $\hz{\Ee}$ to $\hz{\Pp}$.  (See~ \cite[Prop.~8.2 and \S10]{BC92}.)

\begin{lemma}
Let $x_1, \ldots, x_m$ be a regular sequence in a commutative ring~$R$,
and suppose that elements $a_1, \ldots, a_m$ of~$R$ satisfy
\begin{equation}
\sum_{i=1}^m a_i x_i = 0 \quad .
\end{equation}
Then each~$a_i$ lies in the ideal generated by $x_1, \ldots, \hat{x}_i,
\ldots, x_m$.
\end{lemma}

\begin{proof}
We note that the case $m=1$ is trivial, and proceed by induction on~$m$.
The product~$a_m x_m$ lies in the ideal generated by $x_1, \ldots, x_{m-1}$,
and therefore so does~$a_m$ by regularity: say
$a_m = \sum_{i=1}^{m-1} b_i x_i$.  Defining $a'_i = a_i + b_i x_m$
for $1 \leq i \leq m-1$, we have $\sum_{i=1}^{m-1} a'_i x_i = 0$.
Thus we have reduced to the case of~$m-1$, and the result follows
by induction.
\end{proof}

\begin{coroll}
\label{coroll:f_i}
Suppose that $f \in S(\Ee^*)$ lies in the kernel of~$\Inf_{\Ee}^{\Pp}$.
Then there exist elements $f_1, \ldots, f_n$ of~$S(\Ee^*)$ such that
$f = \sum_{i=1}^n f_i \zni{i}$, and the images $\Inf_{\Ee}^{\Pp} (f_i)$
of the~$f_i$ under inflation depend only on~$f$, not on the choice of
the~$f_i$.
\qed
\end{coroll}

\begin{prop}
\label{prop:kappa:p}
There are unique classes $g_0, \ldots, g_{n-1}$ in $\hz{\Pp}$ such that
there exist $f_0, \ldots, f_{n-1}$ in~$S(E^*)$ which satisfy both
$g_i = \Inf_{\Ee}^{\Pp} (f_i)$ and
\begin{equation}
\label{eqn:f_i}
\zni{n+1} + \sum_{i=0}^{n-1} (-1)^{n - i} \zni{i+1} f_i = 0 \quad .
\end{equation}
Each $g_i$~has degree $2 (p^{n+1} - p^{i+1})$,
and $\kappa_{n,r}^p - g_r$ is nilpotent for all $0 \leq r \leq n-1$.
\end{prop}

\begin{proof}
Existence and uniqueness of the~$g_i$ follows from Corollary~\ref{coroll:f_i}.
We once again demonstrate nilpotence by proving that restriction to
every maximal elementary abelian subgroup~$M$ is zero.

Now, each~$\zni{j}$ is a symplectic invariant.  Hence, for any symplectic
transformation $\sigma$~of $\Ee$, Corollary~\ref{coroll:f_i} says that
Eqn.~(\ref{eqn:f_i}) still holds if each~$f_i$ is replaced by~$\sigma^*(f_i)$.
Consequently, $\sigma^*$ fixes each~$g_i$.  As any linear transformation
of a maximal totally isotropic subspace~$I$ may be extended to a symplectic
transformation of~$\Ee$, it follows that each $\Res_M (g_i)$ is a polynomial
in the Dickson invariants for~$I^*$.  Since the symplectic transformations
of~$\Ee$ permute the~$I$ transitively, $\Res_M (g_i)$ is the \emph{same}
polynomial for each~$M$.  Finally, we see by comparing degrees that for
each~$i$ there is a scalar $\lambda_i \in \f$ such that
$\Res_M (g_i) = \lambda_i \DrV{i}{I^*}^p$ for all~$M$.

To establish the result from here, we have to prove that
each $\lambda_i$~is $1$.  First we shall prove that
\begin{equation}
\label{eqn:alternating}
\phi^{p^{n + 1}} + \sum_{i=0}^{n-1} (-1)^{n - i} \phi^{p^{i+1}} g_i = 0
\quad \text{for all $\phi \in \Ee^* \subseteq \hz[2]{\Pp}$.}
\end{equation}
Since the~$g_i$ are invariants for the symplectic group, it suffices to prove
this for one nonzero~$\phi$, say~$\beta_1$.
Recall that the inflation of~$\zni{j}$ is zero, and observe that
differentiating~$\zni{j}$ with respect to~$\alpha_1$ yields~$\beta_1^{p^j}$.
Consequently, differentiating both sides of Eqn.~(\ref{eqn:f_i}) with
respect to~$\alpha_1$ and then inflating yields Eqn.~(\ref{eqn:alternating})
with~$\phi$ replaced by~$\beta_1$.  This establishes
Eqn.~(\ref{eqn:alternating}).

Now restrict to any one maximal elementary abelian $M \leq \Pp$ and
take $p^{\text{th}}$ roots to obtain
\begin{equation}
\phi^{p^n} + \sum_{i=0}^{n-1} (-1)^{n - i} \phi^{p^i}
\lambda_i \DrV{i}{I^*}^p
= 0 \quad \text{for all $\phi \in I^*$.}
\end{equation}
This equation also holds for all $\phi \in I^*$ with each~$\lambda_i$
replaced by~$1$.  Therefore each~$\lambda_i$ must be~$1$, or else taking the
difference of the two left hand sides would yield a polynomial of degree
less than~$p^n$, with too many roots in an integral domain.
\end{proof}

\section{Taking $p^{\text{th}}$ powers of the new classes}
There are $h_1, \ldots, h_{n-1}$ in~$S(\En{n-1}^*)$ such that
$\zni[n-1]{n} = \sum_{i=1}^{n-1} h_i \zni[n-1]{i}$,
by Corollary~\ref{coroll:f_i}.
Now, $S(\En{n-1}^*)$ embeds in $S(\En{n}^*)$, and
$\zni{i} = \zni[n-1]{i} + \zdi{i}$.  Here~$\zdi{i}$ signifies~$\zni[1]{i}$
as a function of $\alpha_n, \beta_n$, that is
$\alpha_n \beta_n^{p^i} - \alpha_n^{p^i} \beta_n$.
So we have
\begin{equation}
\label{eqn:pre_eta}
\zni{n} = \zdi{n} - \sum_{i=1}^{n-1} h_i \zdi{i}
 + \sum_{i=1}^{n-1} h_i \zni{i} \quad .
\end{equation}
As~$\zdi{i}$ is divisible by~$\beta_n$, there is a unique $\eta \in S(\Ee^*)$
such that
\begin{equation}
\label{eqn:eta-prop}
\zni{n} = \beta_n \eta + \sum_{i=1}^{n-1} h_i \zni{i} \quad .
\end{equation}
We may consider~$\eta$ to be an element of $\hz[2p^n]{\Ee}$.
\begin{lemma}
\label{lemma:bottomChi}
For every $\phi \in \proj \Ee^*$, the cohomology class~$\chi_{n-1,\phi}^p$
lies in the image of~$\Inf^{\Pp}_{\Ee}$ modulo nilpotent elements.
\end{lemma}

\begin{proof}
As the symplectic group permutes~$\proj \Ee^*$ transitively, it is enough to
prove the lemma for one~$\phi$.  We shall show that
$\chi_{n-1,\beta_n}^p - \Inf^{\Pp}_{\Ee} (\eta)$ is nilpotent.
Write $\chi$~for $\chi_{n-1,\beta_n}$.

Once more, we prove that the restriction of
$\chi^p - \Inf^{\Pp}_{\Ee} (\eta)$ to every maximal
elementary abelian $M \leq \Pp$ is zero.
The restriction of~$\chi^p$ to each~$M$ is known by
Proposition~\ref{prop:ResChi}.
If the image of~$\beta_n$ in~$I^*$ is nonzero, then it follows from
Eqn.~(\ref{eqn:eta-prop}) that $\Res_M \Inf^{\Pp} (\eta)$ is zero.
Henceforth we may assume that~$M$ lies in the kernel of~$\beta_n$, which
is to say that~$M$ contains~$A_n$ and lies in $\En{n-1} \times \langle
A_n \rangle$.
Then $\Res_M \chi^p = - \Mui{\Ibeta^*}{\alpha_n}^p$ by
Proposition~\ref{prop:ResChi}.

Differentiating Eqn.~(\ref{eqn:eta-prop}) with respect to~$\beta_n$,
then inflating to~$\Pp$ and finally restricting to~$M$, we have
\begin{equation}
- \alpha_n^{p^n} = \Res_M \Inf^{\Pp} (\eta)
- \sum_{i=1}^{n-1} \Res_M \Inf^{\Pp} (h_i) \alpha_n^{p^i} \quad .
\end{equation}
By Proposition~\ref{prop:kappa:p} (for~$\Pn{n-1}$) and
Eqn.~(\ref{eqn:DicksonTheorem}), the choice of the~$h_i$ implies that
\begin{equation}
\alpha_n^{p^n}- \sum_{i=1}^{n-1} \Res_M \Inf^{\Ee} (h_i) \alpha_n^{p^i}
= \Mui{\Ibeta^*}{\alpha_n}^p \quad ,
\end{equation}
whence the result follows.
\end{proof}

\begin{theorem}
\label{thm:pthPower}
For all $0 \leq r \leq n-1$ and all $\phi \in \proj \Ee^*$, the
class~$\chi_{r,\phi}^p$ lies in $\Bild(\Inf^{\Pp}_{\Ee}) + \sqrt{0}$.
\end{theorem}

\begin{proof}
By a result of Evens~\cite{Evens:Steenrod},
Steenrod operations commute with corestriction.  Therefore, at least in
mod-$p$ cohomology, taking $p^{\text{th}}$ powers also commutes with
corestriction.  Hence Lemma~\ref{lemma:bottomChi} and
Proposition~\ref{prop:kappa:p} (applied to~$\Pn{n-1}$) yield the result.
\end{proof}

\section{Linear independence of corestrictions}
We have been looking at three subrings of the cohomology of the extraspecial
$p$-groups.  Firstly the subring generated by top Chern classes, then the
Chern subring, and then the subring generated by top Chern classes and the
corestrictions~$\chi_{r,\phi}$.  In~\cite{p5}, it is shown that top Chern
classes generate a strictly smaller subring than the Chern subring, even
modulo nilradical.  Theorem~\ref{thm:5.2} demonstrates that (at least
modulo nilradical) the Chern subring is contained in the subring generated
by top Chern classes and the~$\chi_{r,\phi}$.
This containment is now shown to be strict.

\begin{prop}
\label{prop:bigger}
Let~$\Ee''$ be a nondegenerate, codimension~$2$ subspace of~$\Ee$.  Exactly
$p+1$ elements $\phi \in \proj \Ee^*$ satisfy $\Ee'' \subseteq
\ker (\phi)$.  When $n=2$, no non-trivial $\f$-linear combination of
these~$\chi_{n-1,\phi}$ lies in $\Bild(\Inf^{\Pp}_{\Ee}) + \sqrt{0}$.
\end{prop}

\begin{coroll}
\label{cor:bigger}
Assume $n=2$.  Modulo nilradical, the Chern subring is strictly contained
in the subring generated by top Chern classes and the~$\chi_{r,\phi}$.
\end{coroll}

\begin{proof}[Proof of Corollary]
Modulo nilradical, the degree~$2(p^2 - p)$ part of the Chern subring consists
of the image of inflation together with one extra class,~$\kappa_1$.
\end{proof}

Let~$\Ee'$ be the orthogonal complement of~$\Ee''$.  Pick a symplectic basis
for~$\Ee$ such that~$\Ee'$ has basis $A_1,B_1$ and $\Ee''$~has basis $A_2,B_2$.
Note that each maximal isotropic subspace $I$~of $\Ee$ satisfies either
$I = (I \cap \Ee') \oplus (I \cap \Ee'')$ or
$(I \cap \Ee') = (I \cap \Ee'') = 0$.
\begin{lemma}
\label{lemma:gamma2}
Consider the $x \in S(\Ee^*)$ whose restriction to~$S(I^*)$ is zero for every
maximal totally isotropic $I \subseteq \Ee$ satisfying
$(I \cap \Ee') = (I \cap \Ee'') = 0$.
These~$x$ form the ideal generated by $\zni[2]{1}$~and
$\gamma_2 = \DrV{1}{{\Ee'}^*} - \DrV{1}{{\Ee''}^*}$.
\end{lemma}

\begin{proof}
Recall that $\zni[2]{1},\zni[2]{2}$ form a regular sequence in~$S(\Ee^*)$,
and generate the joint kernel of restriction to every~$S(I^*)$.
Observe that $(\alpha_1 \beta_1^p - \alpha_1^p \beta_1) \gamma_2
\in \zni[2]{2} + (\zni[2]{1})$.  Moreover, $\gamma_2$~has zero restriction
to~$S(I^*)$ if and only if $I$~satisfies $(I \cap \Ee') = (I \cap \Ee'') = 0$,
whereas $(\alpha_1 \beta_1^p - \alpha_1^p \beta_1)$ restricts to zero
if and only if $I$~satisfies $I = (I \cap \Ee') \oplus (I \cap \Ee'')$.
So the ideal in question is the ideal of~$x$ such that
$(\alpha_1 \beta_1^p - \alpha_1^p \beta_1) x$ has zero restriction to
every~$S(I^*)$.
\end{proof}

\begin{proof}[Proof of Proposition~\ref{prop:bigger}]
Let $\phi \in \proj \Ee^*$ have $\Ee'' \subseteq \ker(\phi)$,
and let~$H$ be the corresponding maximal subgroup of~$\Pp$.  For any
maximal elementary abelian $M \leq H$, have that $I = M / Z$
satisfies $I = (I \cap \Ee') \oplus (I \cap \Ee'')$.  By
Proposition~\ref{prop:ResChi} it follows that $\Res_M (\chi_{1,\phi}) = 0$
for all such~$\phi$ and all~$M$
satisfying $(I \cap \Ee') = (I \cap \Ee'') = 0$.  Invoking
Lemma~\ref{lemma:gamma2} and comparing degrees, the only element of
$\Bild(\Inf^{\Pp}_{\Ee})$ that comes into consideration is~$\gamma_2$ itself.

For no two of these~$\phi$ do the corresponding maximal subgroups~$H$ share
a common maximal elementary abelian subgroup.  Applying
Proposition~\ref{prop:ResChi}
again, these~$\chi_{1,\phi}$ are therefore $\f$-linearly independent.  Moreover,
no nonzero restriction of a~$\chi_{1,\phi}$ is a scalar multiple of the
restriction of~$\gamma_2$.
\end{proof}

\begin{remark}
Proposition~\ref{prop:bigger} and Corollary~\ref{cor:bigger} can be proved
in the same way for general~$n \geq 2$, using the results of~\cite{p5}.
\end{remark}

\section{Chern classes and transfer not sufficient}
Calculating cohomology rings of $p$-groups would be easier if there were a list
of constructions, such as transfer and Chern classes, which together always
yielded a set of generators for the cohomology ring.  Being able to construct
the cohomology ring modulo nilradical would be an important first step.
This gave rise to the so-called $\Chbarop$-conjecture, related to a
construction of Moselle~\cite{Moselle}.  We provide a counterexample to that
conjecture.

\begin{definition}
For a finite group~$G$ and a prime~$p$, the subring $\Chbar{G}$~of $\hf{G}$
is defined recursively as the subring generated by Chern classes for~$G$,
together with the images under corestriction of~$\Chbar{H}$ for proper
subgroups $H$~of $G$.
\end{definition}

\begin{remark}
For $p=2$, one should use Stiefel--Whitney rather than Chern classes to get the
largest possible subring.
\end{remark}

\begin{remark}
Inclusion of $\Chbar{G}$~in $\hf{G}$ always induces an isomorphism of varieties:
see~\cite{VCh}.
\end{remark}

\begin{prop}
\label{prop:Chbar}
Let~$G$ be the extraspecial $3$-group of order~$27$ and exponent~$3$.
Then $\Chbar{G} / \sqrt{0}$ is strictly contained in $\hf{G}/\sqrt{0}$.
\end{prop}

\begin{proof}
Take a symplectic basis $A,B$ for $\Ee = \En{1}$.  We may view~$A,B$ as elements
of~$\Pp = \Pn{1}$.  We shall consider restrictions to the four cyclic subgroups
generated by~$A$, $AB$, $AB^2$, $B$~respectively.  For each of these four
cyclic groups~$K$, denote by~$\xi$ the class in $\hf[1]{K} \cong \Hom(K, \f)$
taking the generator to~$1$.  Let $x \in \mathrm{H}^2$ be the Bockstein of~$\xi$.

For $p$-groups $H_1 < H_2$, recall that $\Cor_{H_1}^{H_2} \Res_{H_1}^{H_2}$ is
zero in mod-$p$ cohomology.
In particular, $\Cor_{H_1}^{H_2}$ is zero whenever~$H_2$ is elementary abelian.
Using the Mackey formula, it follows that $\Res_K^{\Pp} \Cor_H^{\Pp}$ is zero
for any $H < \Pp$ and for $K$~any of the four cyclic subgroups.

Let $\alpha,\beta$ be the dual basis for~$\Ee^*$.  These may be viewed as
elements of $\hf[1]{\Pp}$.  Let $a,b \in \mathrm{H}^2$ be their Bocksteins.
Then $a,b$ are also the first Chern classes of the corresponding
representations, and all first Chern classes are in their span.
Note that the restrictions of~$a$ to the four cyclics are $x,x,x,0$
respectively, whereas the restrictions of~$b$ are $0,x,-x,x$ respectively.

We refer the reader to~\cite{Leary:mod-p} for the fact that $\alpha \beta = 0$
in $\hf[2]{\Pp}$, and for an account of Massey triple products in this context.
Define~$Y$ to be the Massey product
$\langle \alpha, \alpha,\beta\rangle \in \hf[2]{\Pp}$.
By naturality, this has zero restriction
to the cyclics generated by~$A,B$, whereas its restriction
to the subgroup generated by~$AB$ is $\langle\xi,\xi,\xi\rangle$.  This
is~$-x$.  Hence~$Y$ lies outside~$\Chbar{\Pp}$, even modulo nilradical.
\end{proof}

\begin{remark}
We conclude from the proof of Proposition~\ref{prop:Chbar} that any
list of constructions that was always guaranteed to provide enough generators
modulo nilradical for mod-$p$ cohomology would have to include Massey products.
\end{remark}

\begin{remark}
Strictly, Proposition~\ref{prop:Chbar} is only a counterexample to the
mod-$p$ version of the $\Chbarop$-conjecture.  One may however define in the
same way a subring~$\Chbar{G}$ of $\hz{G}$, and it is an open question whether
containment modulo nilradical is strict in this case.
\end{remark}



\begin{thebibliography}{MMM00}







\bibitem[Ben93]{Benson:PolyInvts}
D.~J. Benson.
\newblock \emph{Polynomial Invariants of Finite Groups}.
\newblock London Math.\@ Soc.\@ Lecture Note Ser.\@ no.~190
(Cambridge Univ.\@ Press, 1993).



\bibitem[BC92]{BC92}
D.~J. Benson and J.~F. Carlson.
\newblock The cohomology of extraspecial groups.
\newblock \emph{Bull.\@ London Math.\@ Soc.\@} \textbf{24} (1992), 209--235.
\newblock Erratum: \emph{Bull.\@ London Math.\@ Soc.\@}
\textbf{25} (1993), 498.







\bibitem[Eve68]{Evens:Steenrod}
L.~Evens.
\newblock Steenrod operations and transfer.
\newblock \emph{Proc.\@ Amer.\@ Math.\@ Soc.\@}
\textbf{19} (1968), 1387--1388.

\bibitem[Eve91]{Evens:book}
L.~Evens.
\newblock \emph{The Cohomology of Groups}.
\newblock (Oxford Univ.\@ Press, 1991).

\bibitem[Gre96]{p5}
D.~J. Green.
\newblock ``Chern classes and extraspecial groups of order~$p^5$\@.''
\newblock Submitted for publication.



\bibitem[GL96]{VCh}
D.~J. Green and I.~J. Leary.
\newblock ``The spectrum of the Chern subring.''
\newblock Submitted for publication.





\bibitem[Lea92]{Leary:mod-p}
I.~J. Leary.
\newblock The mod-$p$ cohomology rings of some $p$-groups.
\newblock \emph{Math.\@ Proc.\@ Cambridge Philos.\@ Soc.\@}
\textbf{112} (1992), 63--75.





\bibitem[Mos89]{Moselle}
B.~Moselle.
\newblock ``Calculations in the cohomology of finite groups.''
\newblock Unpublished essay.  Univ.\@ of Cambridge, 1989.

\bibitem[Min95]{Minh:transfer}
P.~A. Minh.
\newblock Transfer map and Hochschild--Serre spectral sequences.
\newblock \emph{J.\@ Pure Appl.\@ Algebra}
\textbf{104} (1995), 89--95.

\bibitem[Mui75]{Mui75}
H.~Mui.
\newblock Modular invariant theory and cohomology
algebras of symmetric groups.
\newblock \emph{J.\@ Fac.\@ Sci.\@ Univ.\@ Tokyo Sect.\@ 1A Math.\@}
\textbf{22} (1975), 316--369.











\bibitem[Tho86]{CBT}
C.~B. Thomas.
\newblock \emph{Characteristic Classes and the Cohomology of Finite Groups}.
\newblock (Cambridge Univ.\@ Press, 1986).




\end{thebibliography}
\end{document}